\documentclass{article}
\usepackage{amssymb,latexsym,amsmath,amsthm,mathrsfs}
\DeclareMathOperator{\Cos}{cos.}
\DeclareMathOperator{\Sin}{sin.}
\DeclareMathOperator{\Sec}{sec.}
\DeclareMathOperator{\Tag}{tag.}
\DeclareMathOperator{\Cot}{cot.}
\begin{document}
\title{Various observations on angles proceeding in geometric
progression\footnote{Presented to the St. Petersburg Academy on November 15, 1773.
Originally published as
{\em Variae observationes circa angulos in progressione geometrica
progredientes},
Opuscula analytica \textbf{1} (1783), 345--352.
E561 in the Enestr{\"o}m index.
Translated from the Latin by Jordan Bell,
Department of Mathematics, University of Toronto, Toronto, Canada.
Email: jordan.bell@gmail.com}}
\author{Leonhard Euler}
\date{}
\maketitle

1. Since a great many well known properties found about angles, or arcs,
and the sines, cosines, tangents, cotangents, secants and cosecants of them
have been derived from the consideration of arcs increasing in arithmetic
progression, those properties which
one may deduce from the consideration of arcs proceeding in geometric
progression seem no less noteworthy, since for the most part the truth of
these seem much more concealed;
therefore I am determined here to unfold many properties of this kind.

2. The starting point that presents itself to us for these types of
speculations is the very well known formula\footnote{Translator: cf.
Thomas Heath,
{\em A history of Greek mathematics}, vol. II, pp. 278--280,
Oxford, 1921.}
\[
\Sin 2\varphi=2\Sin\varphi\cdot\Cos\varphi.
\]
If $s$ denotes any arc or angle, then
\[
\Sin s=2\Sin\tfrac{1}{2}s\cdot \Cos\tfrac{1}{2}s,
\]
and likewise
\[
\Sin\tfrac{1}{2}s=2\Sin\tfrac{1}{4}s\cdot \Cos\tfrac{1}{4}s,
\]
and substituting this value into the former yields
\[
\Sin s=4\Sin\tfrac{1}{4}s\cdot \Cos\tfrac{1}{2}s\cdot \Cos\tfrac{1}{4}s.
\]
Next, because also
\[
\Sin\tfrac{1}{4}s=2\Sin\tfrac{1}{8}s\cdot \Cos\tfrac{1}{8}s,
\]
with this value substituted
\[
\Sin s=8\Sin\tfrac{1}{8}s\cdot \Cos\tfrac{1}{2}s\cdot \Cos\tfrac{1}{4}s\cdot
\Cos\tfrac{1}{8}s.
\]
Proceeding in the same way it is
\[
\Sin s=16\Sin\tfrac{1}{16}s\cdot \Cos\tfrac{1}{2}s\cdot
\Cos\tfrac{1}{4}s\cdot \Cos\tfrac{1}{8}s\cdot \Cos\tfrac{1}{16}s,
\]
and if one continues in this way to infinity, with $i$ denoting an infinite
number or rather an infinite power of 2 we will have
\[
\Sin s=i\Sin\tfrac{s}{i}\cdot \Cos\tfrac{1}{2}s\cdot
\Cos\tfrac{1}{4}s\cdot \Cos\tfrac{1}{8}s\cdot
\Cos\tfrac{1}{16}s \cdot \textrm{etc.},
\]
where, because the arc $\tfrac{s}{i}$ is infinitely small, 
$\Sin\tfrac{s}{i}=\tfrac{s}{i}$, and so 
$i\Sin\tfrac{s}{i}=s$, from which we obtain the outstanding property that
\[
\Sin s=s\Cos\tfrac{1}{2}s\cdot \Cos\tfrac{1}{4}s\cdot \Cos\tfrac{1}{8}s\cdot
\Cos\tfrac{1}{16}s\cdot \Cos\tfrac{1}{32}s\cdot \textrm{etc. to infinity}
\]

3. Therefore the arc $s$ itself can be very prettily defined by its sine
and the cosines of arcs continually diminished in double ratio, as
\[
s=\tfrac{\Sin s}{\Cos \tfrac{1}{2}s \cdot \Cos\tfrac{1}{4}s\cdot
\Cos\tfrac{1}{8}s\cdot \Cos\tfrac{1}{16}s\cdot
\Cos\tfrac{1}{32}s\cdot \textrm{etc.}};
\]
and because $\tfrac{1}{\Cos\varphi}=\Sec\varphi$,
it will be as a complete expression
\[
s=\Sin s\cdot \Sec\tfrac{1}{2}s\cdot
\Sec\tfrac{1}{4}s\cdot \Sec\tfrac{1}{8}s\cdot
\Sec\tfrac{1}{16}s\cdot \Sec\tfrac{1}{32}s\cdot
\Sec\tfrac{1}{64}s\cdot \textrm{etc.},
\]
which expression can be represented quite properly geometrically,
as I have shown before elsewhere.\footnote{Translator: See
Euler's 1738 {\em De
variis modis circuli quadraturam numeris proxime exprimendi},
E74; Opera omnia I.14, p. 257.}

4. Because the arc $s$ is expressed here as a product, by taking logarithms
we will have
\[
ls=l\Sin s+l\Sec\tfrac{1}{2}s+l\Sec\tfrac{1}{4}s+l\Sec\tfrac{1}{8}s
+l\Sec\tfrac{1}{16}s+l\Sec\tfrac{1}{32}s+\textrm{etc.},
\]
whence if we let $s=\tfrac{\pi}{2}=90^\circ$ it will be
\[
l\tfrac{\pi}{2}=0+l\Sec 45^\circ+l\Sec 22^\circ 30'+l\Sec 11^\circ 15'
+l\Sec 5^\circ 37 \tfrac{1}{2}'+\textrm{etc.},
\]
whence, having done a calculation, it will be\footnote{Translator:
Here $l$ denotes the logarithm with base 10. To determine the approximate
sum of the remaining terms, namely $\log \sec \frac{\pi}{2^{k+1}}$
from $k=9$ to $k=\infty$, one can use the Euler-Maclaurin summation formula
applied to the function $f(x)=\log \sec \frac{\pi}{2^{x+1}}$. However,
the integral of $f(x)$ cannot be expressed in a closed form. I don't
know if there are easier ways to approximate this sum.}
\begin{eqnarray*}
l\Sec 45^\circ&=&0,1505150\\
l\Sec 22^\circ 30'&=&0,0343847\\
l\Sec 11^\circ 15'&=&0,0084261\\
l\Sec 5^\circ 37\tfrac{1}{2}'&=&0,0020963\\
l\Sec 2^\circ 48\tfrac{3}{4}'&=&0,0005234\\
l\Sec 1^\circ 24\tfrac{3}{8}'&=&0,0001308\\
l\Sec 0^\circ 42\tfrac{3}{16}'&=&0,0000327\\
l\Sec 0^\circ 21\tfrac{3}{32}'&=&0,0000082\\
\textrm{other terms}&=&0,0000027\\
l\tfrac{\pi}{2}&=&0,1961199\\
l 2&=&0,3010300\\
l \pi&=&0,4971499
\end{eqnarray*}
and hence one finds rather closely that $\pi=3,1415928$.

5. To deduce here new relations, we will differentiate the
last logarithmic equation, and since
\[
d. l\Sec\varphi=\frac{d \varphi \Sin\varphi}{\Cos\varphi}
=d\varphi \Tag\varphi,
\]
the following equation arises by dividing by $ds$
\[
\frac{1}{s}=\Cot s+\tfrac{1}{2}\Tag\tfrac{1}{2}s+
\tfrac{1}{4}\Tag\tfrac{1}{4}s+
\tfrac{1}{8}\Tag\tfrac{1}{8}s+
\tfrac{1}{16}\Tag\tfrac{1}{16}s+
\textrm{etc.},
\]
which series converges very quickly, as will clearly be seen in the following
example. Let us take $s=90^\circ=\tfrac{\pi}{2}$, whence
it will be
\[
\tfrac{2}{\pi}=\tfrac{1}{2}\Tag 45^\circ+\tfrac{1}{4}\tan 22^\circ 30'
+\tfrac{1}{8}\Tag 11^\circ 15'+\tfrac{1}{16}\tan 5^\circ 37\tfrac{1}{2}'
+\textrm{etc.},
\]
and the values taken from tables give
\begin{eqnarray*}
\tfrac{1}{2}\Tag 45^\circ&=&0,5000000\\
\tfrac{1}{4}\Tag 22^\circ 30'&=&0,1035534\\
\tfrac{1}{8}\Tag 11^\circ 15'&=&0,0248640\\
\tfrac{1}{16}\Tag 5^\circ 37\tfrac{1}{2}'&=&0,0061557\\
\tfrac{1}{32}\Tag 2^\circ 48\tfrac{3}{4}'&=&0,0015352\\
\tfrac{1}{64}\Tag 1^\circ 24\tfrac{3}{8}'&=&0,0003836\\
\textrm{for the remaining}&=&0,0001279\\
\tfrac{2}{\pi}&=&0,6366198,
\end{eqnarray*}
hence\footnote{Translator: $\tfrac{1}{0,3183099}=
3.14159251\ldots$ while $\pi=3.14159265\ldots$.}
\[
\pi=\tfrac{2}{0,6366198}=\tfrac{1}{0,3183099}.
\]

6. If we differentiate the last equation for a second time we will come to
a much more convergent series; for since it is\footnote{Translator:
$\Sin\varphi^2=(\Sin\varphi)^2$.}
\[
d.\Cot\varphi=\frac{-d\varphi}{\Sin \varphi^2}
\qquad
\textrm{and}
\qquad
d.\Tag\varphi=\frac{d\varphi}{\Cos\varphi^2}=d\varphi \Sec\varphi^2,
\]
we will get
\[
-\tfrac{1}{ss}=\tfrac{-1}{\Sin s^2}+\tfrac{1}{4}\Sec\tfrac{1}{2}s^2
+\tfrac{1}{16}\Sec\tfrac{1}{4}s^2
+\tfrac{1}{64}\Sec\tfrac{1}{8}s^2+\textrm{etc.}
\]
or
\[
\tfrac{1}{4}\Sec\tfrac{1}{2}s^2+\tfrac{1}{16}\Sec\tfrac{1}{4}s^2
+\tfrac{1}{64}\Sec\tfrac{1}{8}s^2
+\tfrac{1}{256}\Sec\tfrac{1}{16}s^2
+\textrm{etc.}
=\tfrac{1}{\Sin s^2}
-\tfrac{1}{ss}.
\]

7. Let us apply this reasoning to arcs which decrease
in a triple ratio; to this end let us consider
the formula 
\[
\Sin 3\varphi=4\Sin\varphi\Cos\varphi^2-\Sin\varphi
=\Sin\varphi(3-4\Sin\varphi^2),
\]
which gives
\[
\Sin 3\varphi=3\Sin\varphi\bigg(1-\tfrac{4}{3}\Sin\varphi^2\bigg);
\]
so if $s$ denotes any arc, it will be
\[
\Sin s=3\Sin\tfrac{1}{3}s\bigg(1-\tfrac{4}{3}\Sin\tfrac{s^2}{3}\bigg);
\]
and in a similar way
\[
\Sin\tfrac{1}{3}s=3\Sin\tfrac{s}{9}\bigg(1-\tfrac{4}{3}\Sin\tfrac{s^2}{9}\bigg),
\]
so that now
\[
\Sin s=9\Sin\tfrac{1}{9}s\bigg(1-\tfrac{4}{3}\Sin\tfrac{s^2}{3}\bigg)
\bigg(1-\tfrac{4}{3}\Sin\tfrac{s^2}{9}\bigg).
\]

If such substitutions are continued to infinity, as before 
one will be led at last
to this expression
\[
\Sin s=s\bigg(1-\tfrac{4}{3}\Sin\tfrac{s^2}{3}\bigg)
\bigg(1-\tfrac{4}{3}\Sin\tfrac{s^2}{9}\bigg)
\bigg(1-\tfrac{4}{3}\Sin\tfrac{s^2}{27}\bigg) \cdot \textrm{etc.}
\]

8. This factors are too complicated, and one may resolve them into simpler
ones in the following way; namely, because $\Sin \varphi^2=\tfrac{1}{2}-
\tfrac{1}{2}\Cos 2\varphi$,
the general form $1-\tfrac{4}{3}\Sin\varphi^2$ is reduced to
\[
\frac{1}{3}+\frac{2}{3}\Cos 2\varphi,
\] 
which one can write as
\[
\frac{2\Cos 60^\circ+2\Cos 2\varphi}{3}.
\]
Now, since
\[
\Cos a+\Cos b=2\Cos\frac{a+b}{2}\Cos\frac{a-b}{2},
\]
it will be
\[
\Cos 60^\circ + \Cos 2\varphi=2\Cos(30^\circ +\varphi)\Cos(30^\circ-\varphi),
\]
which formula multiplied by $\frac{2}{3}$ yields
\[
1-\frac{4}{3}\Sin \varphi^2=\frac{4}{3}\Cos(30^\circ+\varphi)\Cos(30^\circ-
\varphi).
\]
Whereby if this reduction is applied to all the factors found above,
we will have the following infinite product
\begin{eqnarray*}
\frac{\Sin s}{s}&=&\tfrac{4}{3}\Cos(30^\circ+\tfrac{s}{3})\Cos(30^\circ
-\tfrac{s}{3})\\
&&\tfrac{4}{3}\Cos(30^\circ+\tfrac{s}{9})\Cos(30^\circ
-\tfrac{s}{9})\\
&&\tfrac{4}{3}\Cos(30^\circ+\tfrac{s}{27})\Cos(30^\circ
-\tfrac{s}{27})\\
&&\textrm{etc.},
\end{eqnarray*}
which is expressed by secants
\begin{eqnarray*}
\frac{s}{\Sin s}&=&\tfrac{3}{4}\Sec(30^\circ+\tfrac{s}{3})\Sec(30^\circ
-\tfrac{s}{3})\\
&&\tfrac{3}{4}\Sec(30^\circ+\tfrac{s}{9})\Sec(30^\circ
-\tfrac{s}{9})\\
&&\tfrac{3}{4}\Sec(30^\circ+\tfrac{s}{27})\Sec(30^\circ
-\tfrac{s}{27})\\
&&\textrm{etc.},
\end{eqnarray*}
whose factors approach unity more closely the more the arc $s$ is diminished.

9. If we now take logarithms and differentiate all the terms, the numerical
factors $\frac{3}{4}$, since they are constants, will completely disappear
from the calculation; and since as we saw before
\[
d.l.\Sec \varphi=d\varphi \Tag\varphi,
\]
the following equation will be obtained
\[
\begin{array}{lllllll}
\frac{1}{s}&=&\Cot s&+\tfrac{1}{3}\Tag(30^\circ+\tfrac{s}{3})&+\tfrac{1}{9}
\Tag(30^\circ+\tfrac{s}{9})&+\tfrac{1}{27}\Tag(30^\circ+\tfrac{s}{27})
&+\textrm{etc.}\\
&&&-\tfrac{1}{3}\Tag(30^\circ-\tfrac{s}{3})&-\tfrac{1}{9}\Tag(30^\circ
-\tfrac{s}{9})&-\tfrac{1}{27}\Tag(30^\circ-\tfrac{s}{27})&-\textrm{etc.}
\end{array}
\]
It will be helpful to probe this by an example. Therefore let $s=\frac{\pi}{2}$,
and it will be
{\scriptsize
\[
\begin{array}{lllllll}
\tfrac{2}{\pi}=&\tfrac{1}{3}\Tag 60^\circ&+\tfrac{1}{9}\Tag 40^\circ
&+\tfrac{1}{27}\Tag (33^\circ 20')&+\tfrac{1}{81}\Tag(31^\circ 6\tfrac{2}{3}')
&+\tfrac{1}{243}\Tag(30^\circ 22 \tfrac{2}{9}')&+\textrm{etc.}\\
&-\tfrac{1}{3}\Tag 0^\circ&-\tfrac{1}{9}\Tag 20^\circ
&-\tfrac{1}{27}\Tag(26^\circ 40')
&-\tfrac{1}{81}\Tag(28^\circ 53\tfrac{1}{3}')
&-\tfrac{1}{243}\Tag(29^\circ 37 \tfrac{7}{9}')&-\textrm{etc.}
\end{array}
\]
}

10. This last series seems all the more noteworthy because scarcely
anyone would have been able to demonstrate its truth unless the same method
were used. 
This series is without doubt much 
of a much higher level of investigation
than those
to which we were led by the expansion of the previous case,\footnote{Translator: \S 5.} which was
\[
\frac{1}{s}=\Cot s+\tfrac{1}{2}\Tag\tfrac{1}{2}s+
\tfrac{1}{4}\Tag\tfrac{1}{4}s+
\tfrac{1}{8}\Tag\tfrac{1}{8}s+
\tfrac{1}{16}\Tag\tfrac{1}{16}s+
\textrm{etc.};
\]
the truth of this follows from the well known formula
\[
2\Cot 2\varphi=\Cot\varphi-\Tag\varphi,
\] 
from which we have
\[
\Tag\varphi=\Cot\varphi-2\Cot 2\varphi.
\]
Then if the appropriate values are substituted in place of all the tangents,
the series is arranged in the following way
\[
\frac{1}{s}=
\Bigg\{
\begin{array}{llllll}
\Cot s&+\tfrac{1}{2}\Cot \tfrac{1}{2}s&+\tfrac{1}{4}\Cos\tfrac{1}{4}s
&+\tfrac{1}{8}\Cos\tfrac{1}{8}s&+\cdots&+\tfrac{1}{i}\Cot\tfrac{s}{i}\\
-\Cot s&-\tfrac{1}{2}\Cot \tfrac{1}{2}s&-\tfrac{1}{4}\Cot \tfrac{1}{4}s
&-\tfrac{1}{8}\Cot\tfrac{1}{8}s&-\cdots,&
\end{array}
\]
where one sees that all the terms eliminate each other up to the final term
$\frac{1}{i}\Cot \frac{s}{i}$,
which one can write in the form
\[
\frac{\Cos\tfrac{s}{i}}{i\Sin\tfrac{s}{i}}
\]
with $i$ denoting an infinite number. Since the arc $\frac{s}{i}$ is now
infinitely small, $\Cos \tfrac{s}{i}=1$ and indeed the sine of this arc
will be equal to $\frac{s}{i}$, whence this final term will be
$=\frac{1}{s}$, which is the same value that was found equal to the
series.

11. And for the present case 
the direct demonstration which is to be given is exhibited from that formula
by which the tangent of thrice an angle is expressed.
For if one puts
\[
\Tag \varphi=t,
\]
since one has\footnote{Translator: This by applying the addition
formula 
for $\tan$ twice. The addition formula for $\tan$ is
$\tan(\alpha+\beta)=\frac{\tan\alpha+\tan\beta}{1-\tan\alpha
\tan\beta}$.}
\[
\Tag 3\varphi=\frac{3t-t^3}{1-3tt},
\]
it will be
\[
\Cot 3\varphi=\frac{1-3tt}{3t-t^3}
\qquad \textrm{and} \qquad
3\Cot 3\varphi=\frac{3-9tt}{t(3-tt)}.
\]
Then one subtracts $\Cot\varphi=\tfrac{1}{t}$ and it will be
\[
3\Cot 3\varphi-\Cot\varphi=\frac{-8tt}{t(3-tt)}=\frac{-8t}{3-tt};
\]
then in place of $t$ let us substitute its value $\frac{\Sin\varphi}{\Cos\varphi}$ and we will have
\[
3\Cot 3\varphi-\Cot\varphi=\frac{-8\Sin\varphi\Cos\varphi}{3\Cos\varphi^2
-\Sin\varphi^2}.
\]
Let us deal with the numerator and the denominator of the fraction
in the following way: Since
\[
\Cos\varphi^2=\tfrac{1}{2}+\tfrac{1}{2}\Cos 2\varphi
\qquad \textrm{and} \qquad \Sin\varphi^2=\tfrac{1}{2}-\tfrac{1}{2}
\Cos 2\varphi,
\]
the denominator will take the form $1+2\Cos 2\varphi$ which can therefore
be written as
\[
2\Cos 60^\circ+2\Cos 2\varphi,
\]
which further reduces, because
\[
\Cos a+\Cos b=2\Cos \tfrac{a+b}{2}\Cos \tfrac{a-b}{2},
\]
to
\[
4\Cos(30^\circ+\varphi)\Cos(30^\circ -\varphi).
\]
The numerator is clearly $-4\Sin 2\varphi$, so that we now have
\[
3\Cot 3\varphi-\Cot\varphi=\frac{-\Sin 2\varphi}{\Cos(30^\circ+\varphi)
\Cos(30^\circ-\varphi)}.
\]
Now, since in general
\[
\Sin 2\varphi=\Sin(a+\varphi)\Cos(a-\varphi)
-\Cos(a+\varphi)\Sin(a-\varphi),
\]
let us take $a=30^\circ$ and we will have the following equation
\[
\begin{split}
3\Cot 3\varphi-\Cot\varphi&=
\frac{-\Sin(30^\circ+\varphi)\Cos(30^\circ-\varphi)
+\Cos(30^\circ+\varphi)\Sin(30^\circ-\varphi)}{\Cos(30^\circ+\varphi)
\Cos(30^\circ-\varphi)}\\
&=-\Tag(30^\circ+\varphi)+\Tag(30^\circ-\varphi),
\end{split}
\]
by which we arrive at this noteworthy equation
\[
\Cot 3\varphi=\tfrac{1}{3}\Cot\varphi-\tfrac{1}{3}\Tag(30^\circ
+\varphi)+\tfrac{1}{3}\Tag(30^\circ-\varphi).
\]

12. Now, by writing $s$ in place of $3\varphi$ for our case we
get at once
\[
\Cot s=\tfrac{1}{3}\Cot\tfrac{s}{3}
-\tfrac{1}{3}\Tag(30^\circ+\tfrac{s}{3})+\tfrac{1}{3}
\Tag(30^\circ-\tfrac{s}{3}).
\] 
Indeed in a similar way it will further be
\[
\tfrac{1}{3}\Cot \tfrac{s}{3}=\tfrac{1}{9}\Cot\tfrac{s}{9}
-\tfrac{1}{9}\Tag(30^\circ+\tfrac{s}{9})+\tfrac{1}{9}
\Tag(30^\circ-\tfrac{s}{9}).
\]
Further in the same way,
\[
\tfrac{1}{9}\Cot \tfrac{s}{9}=\tfrac{1}{27}\Cot\tfrac{s}{27}
-\tfrac{1}{27}\Tag(30^\circ+\tfrac{s}{27})+\tfrac{1}{27}
\Tag(30^\circ-\tfrac{s}{27}),
\]
and if we proceed this way to infinity, we will come finally to
a cotangent of the form
\[
\tfrac{1}{i}\Cot\tfrac{s}{i}=\frac{\Cos\tfrac{s}{i}}{i\Sin\tfrac{s}{i}};
\]
hence our equation will take the form\footnote{Translator:
$\tfrac{1}{i}\Cot\tfrac{s}{i}=\frac{\Cos\tfrac{s}{i}}{i\Sin\tfrac{s}{i}}=
\frac{1+O(i^{-2})}{i\cdot (\tfrac{s}{i}+O(i^{-3}))}=
\frac{1+O(i^{-2})}{s+O(i^{-2})}
=\frac{1}{s}+O(i^{-2})$ as $i \to \infty$.}
\[
\begin{array}{llllll}
\Cot s&=-\tfrac{1}{3}\Tag(30^\circ+\tfrac{s}{3})&-\tfrac{1}{9}
\Tag(30^\circ+\tfrac{s}{9})&-\tfrac{1}{27}\Tag(30^\circ+\tfrac{s}{27})&-\cdots
&-\tfrac{1}{s}\\
&+\tfrac{1}{3}\Tag(30^\circ-\tfrac{s}{3})
&+\tfrac{1}{9}\Tag(30^\circ-\tfrac{s}{9})
&+\tfrac{1}{27}\Tag(30^\circ-\tfrac{s}{27})
&+\cdots,&
\end{array}
\]
from which we deduce our very equation itself that was to be demonstrated
\[
\begin{array}{lllllll}
\frac{1}{s}&=&\Cot s&+\tfrac{1}{3}\Tag(30^\circ+\tfrac{s}{3})&+\tfrac{1}{9}
\Tag(30^\circ+\tfrac{s}{9})&+\tfrac{1}{27}\Tag(30^\circ+\tfrac{s}{27})
&+\textrm{etc.}\\
&&&-\tfrac{1}{3}\Tag(30^\circ-\tfrac{s}{3})&-\tfrac{1}{9}\Tag(30^\circ
-\tfrac{s}{9})&-\tfrac{1}{27}\Tag(30^\circ-\tfrac{s}{27})&-\textrm{etc.}
\end{array}
\]

13. Furthermore, one may similarly exhibit for higher ratios series of this type
in  which the arc $s$ is continually diminished.
For since
\[
\Sin 4\varphi=8\Sin \varphi \Cos(45^\circ+\varphi)
\Cos(45^\circ-\varphi)\Cos\varphi,
\]
it will be in a quadruple ratio
\[
\begin{array}{llllll}
\frac{1}{s}&=\Cot s&+\tfrac{1}{4}\Tag \tfrac{s}{4}&+\tfrac{1}{16}\Tag\tfrac{s}{16}
&+\tfrac{1}{64}\Tag\tfrac{s}{64}
&+\textrm{etc.}\\
&&+\tfrac{1}{4}\Tag(45^\circ+\tfrac{s}{4})&
+\tfrac{1}{16}\Tag(45^\circ+\tfrac{s}{16})
&+\tfrac{1}{64}\Tag(45^\circ+\tfrac{s}{64})
&+\textrm{etc.}\\
&&-\tfrac{1}{4}\Tag(45^\circ-\tfrac{s}{4})
&-\tfrac{1}{16}\Tag(45^\circ-\tfrac{s}{16})
&-\tfrac{1}{64}\Tag(45^\circ-\tfrac{s}{64})
&-\textrm{etc.}
\end{array}
\]
Next since\footnote{Translator: See Euler's 1774 {\em Quomodo sinus et
cosinus angulorum multiplorum per producta exprimi queant}, E562;
Opera omnia I.15, p. 509.}
\[
\Sin 5\varphi=16\Sin \varphi \Cos (18^\circ+\varphi)\Cos(18^\circ-\varphi)
\Cos(54^\circ+\varphi)\Cos(54^\circ-\varphi)
\]
we will find in quintuple ratio
\[
\begin{array}{lllll}
\frac{1}{s}&=\Cot s&+\tfrac{1}{5}\Tag(18^\circ+\tfrac{s}{5})
&+\tfrac{1}{25}\Tag(18^\circ+\tfrac{s}{25})&+\textrm{etc.}\\
&&-\tfrac{1}{5}\Tag(18^\circ-\tfrac{s}{5})&-\tfrac{1}{25}\Tag(18^\circ
-\tfrac{s}{25})&-\textrm{etc.}\\
&&+\tfrac{1}{5}\Tag(54^\circ+\tfrac{s}{5})&+\tfrac{1}{25}\Tag(54^\circ
+\tfrac{s}{25})&+\textrm{etc.}\\
&&-\tfrac{1}{5}\Tag(54^\circ-\tfrac{s}{5})&-\tfrac{1}{25}\Tag(54^\circ
-\tfrac{s}{25})&-\textrm{etc.}
\end{array}
\]
One can proceed further in exactly the same way, but truly the resulting
series would be too muddled than would deserve our attention.

\end{document}